\documentclass[10pt]{amsart}       
\usepackage{amssymb}                  
\usepackage{amsthm}
\usepackage{graphicx}
\usepackage{amscd}
\usepackage{color}
\usepackage{enumerate}
\usepackage{xy}
\usepackage{amsmath}
\xyoption{all}

\theoremstyle{plain}
\newtheorem{theorem}{Theorem}

\theoremstyle{definition}

\newtheorem{example}{Example}
\theoremstyle{remark}

\newtheorem{remark}{Remark}
\def\ens{\ensuremath}                  
\newcommand\bydef{\ens{\stackrel{\text{def}}{=}}}
\newcommand{\N}{\mathbb{N}}

\newcommand{\R}{\mathbb{R}}
\newcommand{\T}{\mathbb{T}}
\newcommand{\Z}{\mathbb{Z}}

\newcommand{\card}{\mathbf{card}}                                
\newcommand{\modone}[1]{(\mathrm{mod}\ 1)}            
\newcommand\uvm[2]{\scalebox{#1}{\ens{#2}}} 

\newcommand\hsp[1]{\mbox{}\hspace{#1mm}} 
\newcommand\vsp[1]{\par \vspace{#1mm}} 
\newcommand\vspm[1]{\par \vspace{-#1mm}} 
\newcommand\dist{\text{\rm dist}}  
\def\bfx{{\text{\bf  x}}}                  
\def\bfy{{\text{\bf  y}}}                  
\def\bsq{\hfill\raisebox{-3mm}{$\blacksquare$}\;}

\newcommand\ppp[1]{\mbox{}}                          

\title{The densest sequence in the unit circle.}

\author[M.\ Boshernitzan]{Michael Boshernitzan}

\address{Department of Mathematics, Rice University, Houston, TX~77005, USA}

\email{michael@rice.edu}

\author[J.\ Chaika]{Jon Chaika}

\address{Department of Mathematics, Rice University, Houston, TX~77005, USA}

\email{jmc5124@rice.edu}

\begin{document}
\begin{abstract} 
We exhibit the densest sequence in the unit circle  $\T=\R/\Z$,
$$x_k=\log_2(2k-1)\,  (\text{mod\,}1),\hsp2  k\geq1,$$
in this short note.

 \end{abstract}

\maketitle
\vspm5                                                                                                                \ppp{Section 1\\ Defs 1}                                                                                                                                                                                                                    
\section{Notation and Results.}
Denote by   $\T=\R/\Z$\,  the unit circle and by  $\rho$\,  the standard metric\,
$\rho(s,t)\bydef\dist(s-t,\Z)$\,  on it.  Let   $\T^\N$\,  stand for  the family of sequences\,
$\bfx=\{x_k\}_{k\in\N}$\,  in  $\T$.

For a sequence  $\bfx=\{x_k\}_{k\in\N}\in \T^\N$\,  denote \vsp1
\begin{itemize}
\item  $\bfx[n]=\{x_k\mid 1\leq k\leq n\}$\,  to be  the set  of the first  $n$  elements  of  $\bfx$; \
        $\card\big(\bfx[n]\big)\leq n, \; n\geq1$. \vsp2
\item    $D_n=D_n(\bfx)=\max\limits_{t\in \T}\,\rho(t, \bfx[n])=\max\limits_{t\in \T}
            \big(\min\limits_{1\leq k\leq n}\, \rho(t, x_k)\big),\; n\geq1$. \vsp1
\item    $d_n=d_n(\bfx)=\min\limits_{1\leq r<s\leq n}\  \rho(x_r,x_s)$.
\end{itemize}
\vsp1
 The inequalities  $D_n(\bfx)\geq\uvm{1.2}{\tfrac1{2n}}$  and  $d_n(\bfx)\leq\uvm{1.2}{\frac1n}$  are immediate,  
 for any sequence\,   $\bfx=\{x_k\}_{k\in\N}\in\T^\N$.  The behavior of the sequence  $D_n(\bfx)$  
 as  $n\to\infty$\,  ``measures"  the speed the sequence  $\bfx$  becomes dense.\vsp2 
 
 Denote\, $\phi_n=\uvm{1.15}{\frac{\ln(n+1)-\ln n}{2 \ln 2}}, \quad n\geq1$. 
 The following are main results of this note.

\begin{theorem}\label{thm:1}
For any sequence   $\bfx\in\T^\N$\,  the  inequality\, $D_n(\bfx)\geq \phi_n$\,
has  infinitely many solutions in\, $n\geq1$.  
\end{theorem}

In particular,
$\limsup\limits_{n\to\infty}\, nD_n(\bfx)\geq \uvm{1.15}{\frac1{2\ln2}}$,
for all $\bfx\in\T^\N$.

\begin{theorem}\label{thm:2}
For any sequence   $\bfx\in\T^\N$\,  the  inequality\,
$d_n(\bfx)\!\leq\phi_{n-1}
$\,
has  infinitely many solutions  \text{in \ $n\geq2$}.
\end{theorem}

In particular,\, 
$\liminf\limits_{n\to\infty}\, nd_n(\bfx)\leq\uvm{1.15}{\frac1{2\ln2}}$, 
for all $\bfx\in\T^\N$.
\begin{example} \label{ex:1}
For the sequence\,  $\bfy=\{y_k\}_{k\in\N},\
   y_k=\log_2(2k-1)\,  (\text{\rm mod\ }1)$,\,  one has:
\[
D_n(\bfy)=\phi_n, \quad  
  d_n(\bfy)=2\phi_{2n-1}=\uvm{1.15}{\tfrac{\ln(2n)-\ln(2n-1)}{\ln 2}}, 
  \quad \text{\rm for all } n\geq1.
\]

In particular, $\lim_{n\to\infty} nD_n(\bfy)=\lim_{n\to\infty} n\,d_n(\bfy)=\uvm{1.15}{\frac1{2\ln2}}$.
\end{example}
\vsp1

Example \ref{ex:1}  shows that  Theorem 1 cannot be improved
and that the inequality in Theorem 2 is close to optimal.
We don't know whether Example \ref{ex:1}
provides the optimal upper estimate for  $d_n$ as well.
\vsp2

{\bf Question}.  May  the inequality  $d_n(\bfx)\!<\phi_{n-1}$\,  in Theorem \ref{thm:2}
be replaced by\,  $d_n(\bfx)\!\leq2\phi_{2n-1}?$
\begin{remark}

These theorems and example deal with the most economical packing for $\T=[0,1)$  
by a sequence of points  $\bfx=\{x_k\}_1^\infty$ in $\T$   which is optimal over all stopping times  $n\in\N$.
 A close example is the packing of sunflowers (\cite[pages 112--114]{livio});  this problem also considers packings  ``at all times". However, it restricts itself to packings given by the orbit of a rotation. Our example provides different results (a better constant,  with $d_n$  strictly decreasing)  because no restrictions are imposed on  $\bfx$.
 Optimal packings in other situations, with other constraints, have also previously been dealt with (see \cite{hales}).

\end{remark}
\vsp3

\section{The proofs.}
\noindent {\bf Proof of Theorem \ref{thm:1}}.
WLOG  we assume that  $x_1=0$  and that  $x_i\neq x_j$  for
$i\neq j$,  i.\,e. that\,  $\card(\bfx[n])=n$.  For  $n\geq1$,  denote  by  $J_{n,1}, J_{n,2}, \ldots, J_{n,n}$  
the  $n$  intervals  corresponding to the partition of\,  $\T=[0,1)$  by the set  $\bfx[n]$.  The intervals
$J_{n,k}$  are assumed to be arranged in the non-increasing order according to their lengths   
$t_{n,k}=|J_{n,k}|$: 
$
t_{n,1}\geq t_{n,2}\geq\ldots\geq t_{n,n}.
$
Observe that   
\begin{equation}\label{eq:ineq}
t_{n,k}\leq t_{n+1,k-1}, \ \text{ for }\,  1<k\leq n,
\end{equation}
and hence
$
t_{n,k}\leq t_{n+k-1,1}.
$
To complete the proof of the theorem,  assume  to the contrary  that\,
   $D_n<\phi_n$\,  for all large   $n$.
Then
$$
1=\sum_{k=1}^n  t_{n,k}\leq\sum_{k=1}^n  t_{n+k-1,1}=\sum_{k=1}^n  2D_{n+k-1}<2\sum_{k=1}^n \phi_{n+k-1}=1,
$$
   a contradiction.       \bsq
\vsp5
   
\noindent {\bf Proof of Theorem \ref{thm:2}}.  Define 
$J_{n,k}$  and  $t_{n,k},\ 1\leq k\leq n$,   just as   in the proof of   
Theorem \ref{thm:1}. (Again, we assume  that   $x_1=0$  and that  $\card(\bfx[n])=n$).

The inequality  \eqref{eq:ineq}  implies that\,
$t_{2n,2k}\geq t_{n+k,n+k}$,\,  for  $1\leq k\leq n.$

Assume to the contrary that  $d_n>\phi_{n-1}$   for all large  $n$.
A contradiction is derived as follows:
\begin{align*}
1/2=\sum_{k=2n}^{4n-1} \phi_k=&\sum_{k=2n+1}^{4n} \phi_{k-1}\hsp{2.5}<
\sum_{k=2n+1}^{4n} d_k=\sum_{k=2n+1}^{4n} t_{k,k}=\\
=&\sum_{k=1}^{2n}t_{2n+k,2n+k}\leq \sum_{k=1}^{2n} t_{4n,2k}\leq
\uvm{1.1}{\tfrac12} \sum_{k=1}^{4n} t_{4n,k}=1/2.
\end{align*}

\bsq
\vsp5

\end{document}